# A Geometric Proof of Generalized Riemann Hypothesis

**Kaida Shi**

Department of Mathematics, Zhejiang Ocean University,

Zhoushan City, 316004, Zhejiang Province, China

**Abstract** Beginning from the resolution of Dirichlet $L$ function $L(s, \chi)$, using the inner product formula of infinite-dimensional vectors in the complex space, the author proved the world's baffling problem generalized Riemann hypothesis.

MR(2000): 11M

## 1 Introduction

For investigating the distribution problem of prime numbers within the arithmetic series, in the first, L. Dirichlet imported the Dirichlet $L$ function $L(s, \chi)$. Although its property and effect are similar with Riemann Zeta function $\zeta(s)$, but the different matter is, the research about the distribution of the real zeros of the Dirichlet $L$ function $L(s, \chi)$ which correspond real property is a very difficult thing. But because of this, it has very important meaning.

For solving above-mentioned problem, L. Dirichlet put forward a famous proposition as follows:

**Generalized Riemann hypothesis (GRH):** For any $\chi \bmod q$, all the non-trivial zeroes of the Dirichlet $L$ function $L(s, \chi)$ $(\operatorname{Re}(s) > 1)$ lie on the straight line $\operatorname{Re}(s) = \frac{1}{2}$.

In the famous paper "Ueber die Anzahl der Primzahlen unter einer gegebenen Große", B. Riemann revised $s$ as $s = \sigma + it$. Also, for Dirichlet $L$ function $L(s, \chi)$ $(\operatorname{Re}(s) > 1)$, we revised $s$ as $s = \sigma + it$.

Now, let's prove the proposition.

## 2 The determination of real part of non-trivial zeroes of the Dirichlet $L$ function $L(s, \chi)$

According to the **symmetry** (about $\operatorname{Re}(s) = \frac{1}{2}$) of distribution of non-trivial zeroes of the Dirichlet $L$

function $L(s, \chi)$, we must consider the **uniqueness** and **continuity** of non-trivial zeroes of the Dirichlet $L$ function $L(s, \chi)$ for both the imaginative coordinate $t_0$ and the real coordinate $\sigma_0$, therefore, we suggest:

**Theorem** *The imaginative coordinate $t_0$ of each non-trivial zero of the Dirichlet L function $L(s, \chi)$ doesn't correspond with two or over two real coordinates $\sigma_k (k = 1, 2, 3, \cdots)$; but the real coordinate $\sigma_0$ of each non-trivial zero of the Dirichlet L function $L(s, \chi)$ corresponds with infinite imaginative coordinates $t_k (k = 1, 2, 3, \cdots)$.*

**Proof** First, suppose that $t_0$ corresponds with two real coordinates $\sigma_1$ and $\sigma_2$ $(\sigma_1 \neq \sigma_2)$, then, we have following equation group:

$$\begin{cases} L(\sigma_1 + it_0, \chi) = \dfrac{\chi(1)}{1^{\sigma_1+it_0}} + \dfrac{\chi(2)}{2^{\sigma_1+it_0}} + \dfrac{\chi(3)}{3^{\sigma_1+it_0}} + \cdots + \dfrac{\chi(n)}{n^{\sigma_1+it_0}} + \cdots = 0, \\ L(\sigma_2 + it_0, \chi) = \dfrac{\chi(1)}{1^{\sigma_2+it_0}} + \dfrac{\chi(2)}{2^{\sigma_2+it_0}} + \dfrac{\chi(3)}{3^{\sigma_2+it_0}} + \cdots + \dfrac{\chi(n)}{n^{\sigma_2+it_0}} + \cdots = 0. \end{cases}$$

Taking the first equation minus the second equation, we will obtain

$$L(\sigma_1 + it_0, \chi) - L(\sigma_2 + it_0, \chi) = \sum_{n=1}^{\infty} (\dfrac{\chi(n)}{n^{\sigma_1+it_0}} - \dfrac{\chi(n)}{n^{\sigma_2+it_0}})$$

$$= \sum_{n=1}^{\infty} \dfrac{n^{\sigma_2} - n^{\sigma_1}}{n^{(\sigma_1+\sigma_2)+it_0}} \chi(n),$$

because $\sigma_1 \neq \sigma_2$, therefore, $n^{\sigma_2} - n^{\sigma_1} \neq 0$. On the other hand, because

$$n^{(\sigma_1+\sigma_2)+it_0} = n^{\sigma_1+\sigma_2} \cdot n^{it_0} = n^{\sigma_1+\sigma_2} \cdot e^{it_0 \ln n}$$
$$= n^{\sigma_1+\sigma_2} \cdot (\cos(t_0 \ln n) + i\sin(t_0 \ln n)) \neq 0,$$

therefore we have

$$L(\sigma_1 + it_0, \chi) - L(\sigma_2 + it_0, \chi) \neq 0$$

namely

$$L(\sigma_1 + it_0, \chi) \neq L(\sigma_2 + it_0, \chi).$$

This is contradictory with the fact that the first equation minus the second equation must equal to zero. Hence, the uniqueness of the real coordinate $\sigma$ of each non-trivial zero of the Dirichlet $L$ function $L(s, \chi)$ has been proved.

Second, suppose that the real coordinate $\sigma_0$ of each non-trivial zero of the Dirichlet $L$ function $L(s, \chi)$ corresponds with two imaginative coordinates $t_1$ and $t_2$, then, we have following equation group:

$$\begin{cases} L(\sigma_0 + it_1, \chi) = \dfrac{\chi(1)}{1^{\sigma_0+it_1}} + \dfrac{\chi(2)}{2^{\sigma_0+it_1}} + \dfrac{\chi(3)}{3^{\sigma_0+it_1}} + \cdots + \dfrac{\chi(n)}{n^{\sigma_0+it_1}} + \cdots = 0, \\ L(\sigma_0 + it_2, \chi) = \dfrac{\chi(1)}{1^{\sigma_0+it_2}} + \dfrac{\chi(2)}{2^{\sigma_0+it_2}} + \dfrac{\chi(3)}{3^{\sigma_0+it_2}} + \cdots + \dfrac{\chi(n)}{n^{\sigma_0+it_2}} + \cdots = 0. \end{cases}$$

Taking the first equation minus the second equation, we obtain

$$L(\sigma_0 + it_1, \chi) - L(\sigma_0 + it_2, \chi) = \sum_{n=1}^{\infty} \left( \dfrac{\chi(n)}{n^{\sigma_0+it_1}} - \dfrac{\chi(n)}{n^{\sigma_0+it_2}} \right)$$

$$= \sum_{n=1}^{\infty} \dfrac{n^{it_2} - n^{it_1}}{n^{\sigma_0+i(t_1+t_2)}} \chi(n) = \sum_{n=1}^{\infty} \dfrac{e^{it_2 \ln n} - e^{it_1 \ln n}}{n^{\sigma_0+i(t_1+t_2)}} \chi(n)$$

$$= \sum_{n=1}^{\infty} \dfrac{(\cos(t_2 \ln n) - \cos(t_1 \ln n)) + i(\sin(t_2 \ln n) - \sin(t_1 \ln n))}{n^{\sigma_0+i(t_1+t_2)}} \chi(n),$$

where

$$n^{\sigma_0+i(t_1+t_2)} = n^{\sigma_0} \cdot n^{i(t_1+t_2)} = n^{\sigma_0} \cdot e^{i(t_1+t_2)\ln n}$$
$$= n^{\sigma_0} \cdot (\cos(t_1 \ln n) + i\sin(t_1 \ln n)) \cdot (\cos(t_2 \ln n) + i\sin(t_2 \ln n))$$
$$\neq 0.$$

Enable above expression equals to zero, we must have

$$\begin{cases} \cos(t_2 \ln n) = \cos(t_1 \ln n), \\ \sin(t_2 \ln n) = \sin(t_1 \ln n) \end{cases} \quad (n = 1, 2, 3, \cdots)$$

so, we obtain $t_1 = t_2 + \dfrac{2k\pi}{\ln n}$ $(k = 1, 2, 3, \cdots)$. That is to say $t_1$ and $t_2$ can take any value, but

$t_1 - t_2 = \dfrac{2k\pi}{\ln n}$ $(k = 1, 2, 3, \cdots)$.

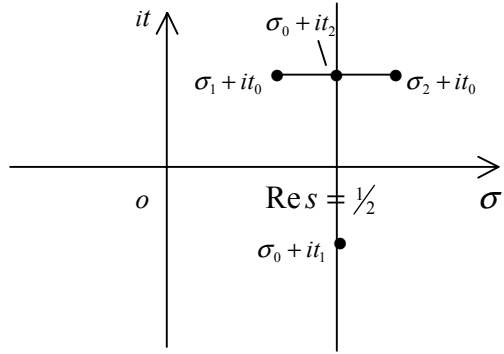

Fig. 1

On the other hand, we have the Dirichlet $L$ functional equality:

$$L(1-s, \chi) = (2\pi)^{-s} q^{s-1} \tau(\chi) \Gamma(s) \cdot \left\{ e^{-i\pi s/2} + \chi(-1) e^{i\pi s/2} \right\} \cdot L(s, \overline{\chi}). \tag{1}$$

where $\chi(-1) = \pm 1$ and $q = \chi(-1)\tau(\chi)\tau(\overline{\chi})$, therefore, the expression (1) becomes:

$$L(1-s, \chi) = 2 \cdot (2\pi)^{-s} q^{s-1} \tau(\chi) \Gamma(s) \cdot \cos\dfrac{\pi s}{2} \cdot L(s, \overline{\chi}). \tag{2}$$

and

$$L(1-s, \chi) = -2i \cdot (2\pi)^{-s} q^{s-1} \tau(\chi) \Gamma(s) \cdot \sin\dfrac{\pi s}{2} \cdot L(s, \overline{\chi}) \tag{3}$$

From (2), if and only if $s = \dfrac{1}{2}$, we have

$$L(\tfrac{1}{2},\chi) = 2 \cdot (2\pi)^{-\tfrac{1}{2}} q^{-\tfrac{1}{2}} \tau(\chi) \Gamma(\tfrac{1}{2}) \cdot \cos\frac{\pi}{4} \cdot L(\tfrac{1}{2},\overline{\chi})$$

$$= \sqrt{2}\,\frac{\tau(\chi)\Gamma(\tfrac{1}{2})}{\sqrt{2\pi q}} L(\tfrac{1}{2},\overline{\chi}) = \frac{\tau(\chi)\Gamma(\tfrac{1}{2})}{\sqrt{\pi\chi(-1)\tau(\chi)\tau(\overline{\chi})}} L(\tfrac{1}{2},\overline{\chi}) =$$

$$= \frac{\sqrt{\tau(\chi)}\Gamma(\tfrac{1}{2})}{\sqrt{\pi\chi(-1)\tau(\overline{\chi})}} L(\tfrac{1}{2},\overline{\chi}).$$

Because $\Gamma(\tfrac{1}{2}) = \sqrt{\pi}$ and take $\chi(-1) = 1$, therefore, we have

$$L(\tfrac{1}{2},\chi) = \frac{\sqrt{\tau(\chi)}}{\sqrt{\tau(\overline{\chi})}} L(\tfrac{1}{2},\overline{\chi}).$$

If and only if $\chi = \overline{\chi}$, we have

$$L(\tfrac{1}{2},\chi) = L(\tfrac{1}{2},\chi).$$

From (3), if and only if $s = \dfrac{1}{2}$, we have

$$L(\tfrac{1}{2},\chi) = -2i \cdot (2\pi)^{-\tfrac{1}{2}} q^{-\tfrac{1}{2}} \tau(\chi) \Gamma(\tfrac{1}{2}) \cdot \sin\frac{\pi}{4} \cdot L(\tfrac{1}{2},\overline{\chi})$$

$$= -\sqrt{2}i \cdot \frac{\tau(\chi)\Gamma(\tfrac{1}{2})}{\sqrt{2\pi q}} L(\tfrac{1}{2},\overline{\chi}) = -i \cdot \frac{\tau(\chi)\Gamma(\tfrac{1}{2})}{\sqrt{\pi\chi(-1)\tau(\chi)\tau(\overline{\chi})}} L(\tfrac{1}{2},\overline{\chi}) =$$

$$= -i \cdot \frac{\sqrt{\tau(\chi)}\Gamma(\tfrac{1}{2})}{\sqrt{\pi\chi(-1)\tau(\overline{\chi})}} L(\tfrac{1}{2},\overline{\chi}).$$

Because $\Gamma(\tfrac{1}{2}) = \sqrt{\pi}$ and take $\chi(-1) = -1$, we have

$$L(\tfrac{1}{2},\chi) = -i \cdot \frac{\sqrt{\tau(\chi)}}{i \cdot \sqrt{\tau(\overline{\chi})}} L(\tfrac{1}{2},\overline{\chi}),$$

therefore, if and only if $\chi = \overline{\chi}$, we have

$$L(\tfrac{1}{2},\chi) = -L(\tfrac{1}{2},\chi).$$

From the equation group

$$\begin{cases} L(\tfrac{1}{2},\chi) = L(\tfrac{1}{2},\chi); \\ L(\tfrac{1}{2},\chi) = -L(\tfrac{1}{2},\chi), \end{cases}$$

we can know

$$L(\tfrac{1}{2},\chi) = 0.$$

**So, we have known the real part** $\mathrm{Re}(s) = \sigma$ **of the non-trivial zeroes of the Dirichlet L function** $L(s,\chi)$ **is** $\dfrac{1}{2}$.

## 3  The resolution of Dirichlet L function $L(s,\chi)$

Considering the geometric meaning of the Dirichlet L function

$$L(s,\chi) = \sum_{n=1}^{\infty} \frac{\chi(n)}{n^s} = \frac{\chi(1)}{1^s} + \frac{\chi(2)}{2^s} + \frac{\chi(3)}{3^s} + \cdots + \frac{\chi(n)}{n^s} + \cdots, \qquad (4)$$

we have the following figure:

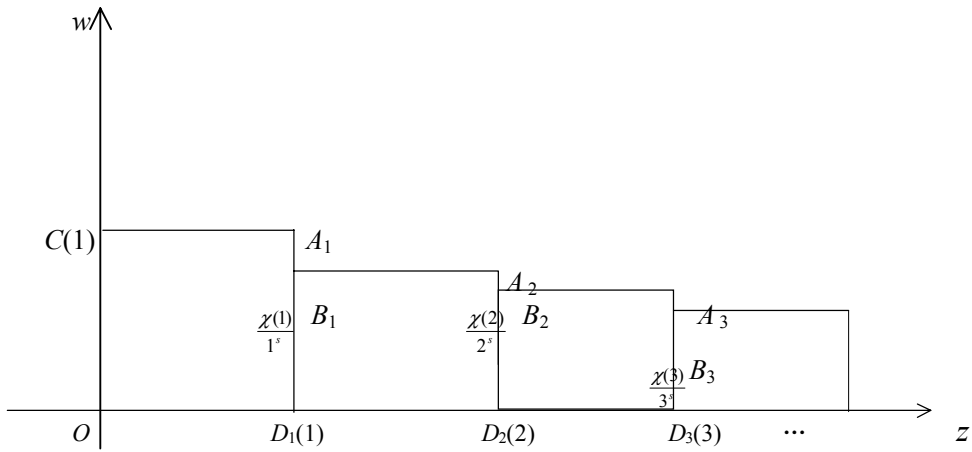

Fig. 2

In this figure, the areas of the rectangles $A_1 D_1 O C, A_2 D_2 D_1 B_1, A_3 D_3 D_2 B_2, \cdots$ are respectively:

$$\frac{\chi(1)}{1^s} \cdot 1, \ \frac{\chi(2)}{2^s} \cdot 1, \ \frac{\chi(3)}{3^s} \cdot 1, \ \cdots, \ \frac{\chi(n)}{n^s} \cdot 1, \ \cdots$$

therefore, the geometric meaning of the Dirichlet $L$ function $L(s,\chi)$ is the sum of the areas of a series of rectangles within the complex space $s$.

Using the inner product formula between two infinite-dimensional vectors, the Dirichlet $L$ function $L(s,\chi)$ equation

$$L(s,\chi) = \frac{\chi(1)}{1^s} + \frac{\chi(2)}{2^s} + \frac{\chi(3)}{3^s} + \cdots + \frac{\chi(n)}{n^s} + \cdots = 0, \quad (s = \sigma + it)$$

can be resolved as

$$L(s,\chi) = (\frac{\chi(1)}{1^\sigma}, \frac{\chi(2)}{2^\sigma}, \frac{\chi(3)}{3^\sigma}, \cdots, \frac{\chi(n)}{n^\sigma}, \cdots) \cdot (\frac{1}{1^{it}}, \frac{1}{2^{it}}, \frac{1}{3^{it}}, \cdots, \frac{1}{n^{it}}, \cdots) = 0. \quad (5)$$

or

$$L(s,\chi) = (\frac{1}{1^\sigma}, \frac{1}{2^\sigma}, \frac{1}{3^\sigma}, \cdots, \frac{1}{n^\sigma}, \cdots) \cdot (\frac{\chi(1)}{1^{it}}, \frac{\chi(2)}{2^{it}}, \frac{\chi(3)}{3^{it}}, \cdots, \frac{\chi(n)}{n^{it}}, \cdots) = 0. \quad (6)$$

From the expression (5), we obtain:

$$\sqrt{\sum_{n=1}^{\infty} \frac{\chi^2(n)}{n^{2\sigma}}} \cdot \sqrt{\sum_{n=1}^{\infty} \frac{1}{n^{2it}}} \cdot \cos(\widehat{\vec{N}_1, \vec{N}_2}) = 0. \quad (7)$$

But in the complex space, if the inner product between two vectors equals to zero, then these two vectors are perpendicular, namely, $(\widehat{\vec{N}_1, \vec{N}_2}) = \frac{\pi}{2}$, therefore

$$\cos(\widehat{\vec{N}_1, \vec{N}_2}) = 0.$$

So long as the values of two radical expressions in the left side of (7) are finite (real or imaginary), then we think the expression (7) is tenable.

## 4 The relationship between the volume of the rotation solid and the area of its axis-cross section within the complex space

From the expression (7), we can know that when $\mathrm{Re}(s) = \sigma = \frac{1}{2}$, the series within the first radical expression of left side is **probably** divergent. Because $\cos(\widehat{\vec{N}_1, \vec{N}_2}) = 0$, therefore the situation of the complex series within the second radical expression of left side will change unable to research. Hence, we must transform the Dirichlet $L$ function $L(s, \chi)$ equation. For this aim, let's derive the relationship between the volume of the rotation solid and the area of its axis-cross section within the complex space.

We call the cross section which pass through the axis $z$ and intersects the rotation solid as **axis-cross section**. According to the barycentric formula of the complex plane lamina:

$$\begin{cases} \xi = \dfrac{\int_a^b z[f(z)-g(z)]dz}{\int_a^b [f(z)-g(z)]dz}, \\ \eta = \dfrac{\frac{1}{2}\int_a^b [f^2(z)-g^2(z)]dz}{\int_a^b [f(z)-g(z)]dz}. \end{cases}$$

we have

$$2\pi\eta = \dfrac{\pi\int_a^b [f^2(z)-g^2(z)]dz}{\int_a^b [f(z)-g(z)]dz}. \qquad (*)$$

The numerator of the fraction of right side of the formula (*) is the volume of the rotation solid, and the denominator is the area of the axis-cross section. The $\eta$ of left side of the formula (*) is the longitudinal coordinate of the barycenter of axis-cross section. Its geometric explanation is: taking the axis-cross section around the axis $z$ to rotate the angle of $2\pi$, we will obtain the volume of the rotation solid.

## 5 The proof of divergence of two concerned serieses

Now, let's prove the divergence of the series $\sum\limits_{n=1}^{\infty} \cos(4t \ln n)$ and the series $\sum\limits_{n=1}^{\infty} \sin(4t \ln n)$.

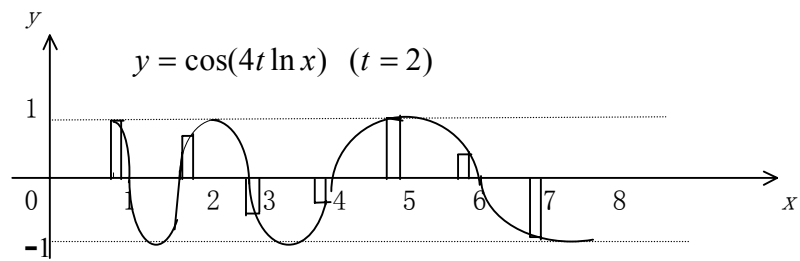

Fig. 3

Obviously, from Fig.3, we can find that the area of $c \cdot \sum_{n=1}^{\infty} \cos(4t \ln n)$ is contained by the area of

$$\int_{1}^{+\infty} \cos(4t \ln x) dx.$$

Suppose that $c$ is a very small positive number, then the sum of the areas of a series of rectangles (which take $\cos(4t \ln n)(n = 1,2,3,\cdots)$ as **high** and take $c$ as **width**) is:

$$c \cdot \sum_{n=1}^{\infty} \cos(4t \ln n) = c \cdot \cos(4t \ln 1) + c \cdot \cos(4t \ln 2) + c \cdot \cos(4t \ln 3) + \cdots.$$

Because

$$\int_{1}^{+\infty} \cos(4t \ln x) dx = \lim_{\substack{w \to +\infty \\ t \to \infty}} \frac{1}{1+16t^2} x[\cos(4t \ln x) + 4t \sin(4t \ln x)]_{1}^{w}$$

$$= \lim_{\substack{w \to +\infty \\ t \to \infty}} \frac{1}{1+16t^2} [w(\cos(4t \ln w) + 4t \sin(4t \ln w)) - 1]$$

$$= \lim_{\substack{w \to +\infty \\ t \to \infty}} \frac{w}{\sqrt{1+16t^2}} \left( \frac{1}{\sqrt{1+16t^2}} \cos(4t \ln w) + \frac{4t}{\sqrt{1+16t^2}} \sin(4t \ln w) \right) - \lim_{t \to \infty} \frac{1}{1+16t^2}$$

Suppose that $\dfrac{1}{\sqrt{1+16t^2}} = \sin \delta$ and $\dfrac{4t}{\sqrt{1+16t^2}} = \cos \delta$, then above equals to

$$\lim_{\substack{w \to +\infty \\ t \to \infty}} \frac{w}{\sqrt{1+16t^2}} \sin(\delta + 4t \ln w) - \lim_{t \to \infty} \frac{1}{1+16t^2}$$

$$= \lim_{\substack{w \to +\infty \\ t \to \infty}} \frac{d}{\sqrt{\frac{1}{t^2}+16}} \sin(\delta + 4t \ln w)$$

$$= \frac{d}{4} \lim_{\substack{w \to +\infty \\ t \to \infty}} \sin(\delta + 4t \ln w).$$

Because the relationship between $w$ and $t$ is linear, therefore, we denote $w = dt$.

Obviously,

$$-\frac{d}{4} \leq \frac{d}{4} \lim_{\substack{w \to +\infty \\ t \to \infty}} \sin(\delta + 4t \ln w) \leq \frac{d}{4},$$

therefore

$$-\frac{d}{4} \leq c \cdot \sum_{n=1}^{\infty} \cos(4t \ln n) \leq \frac{d}{4},$$

namely,

$$-\frac{d}{4c} \leq \sum_{n=1}^{\infty} \cos(4t \ln n) \leq \frac{d}{4c}.$$

Thus, although the series $\sum_{n=1}^{\infty} \cos(4t \ln n)$ is divergent, but it is **oscillatory** and **bounded**.

Similarly, although the series $\sum_{n=1}^{\infty} \sin(4t \ln n)$ is divergent, but it is also **oscillatory** and **bounded**.

In a word, when $t = 0$, the series $\sum_{n=1}^{\infty} \cos(4t \ln n)$ is divergent (goes to **infinite**); when $t \neq 0$, it is **oscillatory** and **bounded**.

Thus, the divergence (**oscillatory** and **bounded**) of $\sum_{n=1}^{\infty} \cos(4t \ln n)$ has been proved.

Similarly, we can prove the divergence (**oscillatory** and **bounded**) of the series $\sum_{n=1}^{\infty} \sin(4t \ln n) \ (t \neq 0)$.

## 6 The proof of generalized Riemann hypothesis

Because we have the volume formula of the rotation solid formed by rotating the rectangle $A_n D_n D_{n-1} B_{n-1}$ around the axis $z$:

$$V_n = \int_{n-1}^{n} \frac{\pi \chi^2(n)}{n^{2s}} dz = \frac{\pi \chi^2(n)}{n^{2s}} z \Big|_{n-1}^{n} = \frac{\pi \chi^2(n)}{n^{2s}},$$

therefore, the sum of the volumes of a series of the cylinders formed by rotating a series of the rectangles $A_1D_1OC, A_2D_2D_1B_1, A_3D_3D_2B_2, \cdots$ around the axis $z$ is:

$$V = \sum_{n=1}^{\infty} V_n = \pi(\frac{\chi^2(1)}{1^{2s}} + \frac{\chi^2(2)}{2^{2s}} + \frac{\chi^2(3)}{3^{2s}} + \cdots + \frac{\chi^2(n)}{n^{2s}} + \cdots).$$

But the sum of the areas of a series of the rectangles $A_1D_1OC, A_2D_2D_1B_1, A_3D_3D_2B_2, \cdots$ is:

$$\frac{\chi(1)}{1^s} + \frac{\chi(2)}{2^s} + \frac{\chi(3)}{3^s} + \cdots + \frac{\chi(n)}{n^s} + \cdots.$$

By the formula (*), we have

$$\pi(\frac{\chi^2(1)}{1^{2s}} + \frac{\chi^2(2)}{2^{2s}} + \frac{\chi^2(3)}{3^{2s}} + \cdots + \frac{\chi^2(n)}{n^{2s}} + \cdots)$$

$$= 2\pi\eta(\frac{\chi(1)}{1^s} + \frac{\chi(2)}{2^s} + \frac{\chi(3)}{3^s} + \cdots + \frac{\chi(n)}{n^s} + \cdots). \tag{8}$$

Substituting the Dirichlet $L$ function $L(s, \chi)$ equation

$$L(s, \chi) = \frac{\chi(1)}{1^s} + \frac{\chi(2)}{2^s} + \frac{\chi(3)}{3^s} + \cdots + \frac{\chi(n)}{n^s} + \cdots = 0$$

into (8), we can obtain the transformed equation:

$$\frac{\chi^2(1)}{1^{2s}} + \frac{\chi^2(2)}{2^{2s}} + \frac{\chi^2(3)}{3^{2s}} + \cdots + \frac{\chi^2(n)}{n^{2s}} + \cdots = 0.$$

This can be explained geometrically as: **if the area of the axis-cross section equals to zero, then it corresponds to the volume of the rotation solid also equals to zero**.

Because $s = \sigma + it$, therefore

$$\frac{\chi^2(1)}{1^{2(\sigma+it)}} + \frac{\chi^2(2)}{2^{2(\sigma+it)}} + \frac{\chi^2(3)}{3^{2(\sigma+it)}} + \cdots + \frac{\chi^2(n)}{n^{2(\sigma+it)}} + \cdots = 0.$$

According to the inner product formula between two infinite-dimensional vectors, above expression can be written as:

$$(\frac{\chi^2(1)}{1^{2\sigma}}, \frac{\chi^2(2)}{2^{2\sigma}}, \frac{\chi^2(3)}{3^{2\sigma}}, \cdots, \frac{\chi^2(n)}{n^{2\sigma}}, \cdots) \cdot (\frac{1}{1^{2it}}, \frac{1}{2^{2it}}, \frac{1}{3^{2it}}, \cdots, \frac{1}{n^{2it}}, \cdots) = 0.$$

Denoting above two infinite-dimensional vectors as $\vec{N}_3$ and $\vec{N}_4$ respectively, from the expression (5), we obtain:

$$\sqrt{\sum_{n=1}^{\infty} \frac{\chi^4(n)}{n^{4\sigma}}} \cdot \sqrt{\sum_{n=1}^{\infty} \frac{1}{n^{4it}}} \cdot \cos(\widehat{\vec{N}_3, \vec{N}_4}) = 0. \qquad (9)$$

Substituting it into the expression (9), we obtain:

$$\sqrt{\sum_{n=1}^{\infty} \frac{\chi^4(n)}{n^2}} \cdot \sqrt{\sum_{n=1}^{\infty} \frac{1}{n^{4it}}} \cdot \cos(\widehat{\vec{N}_3, \vec{N}_4}) = 0.$$

According to the definition of the character $\chi(n)$ of the module $q$, when $(n, q) = 1$, $|\chi(n)| = 1$. Suppose that $\chi(n) = \alpha_n + i\beta_n$, where $\alpha_n^2 + \beta_n^2 = 1$, $-1 \leq \alpha_n \leq 1$, $-1 \leq \beta_n \leq 1$. We have

$$\chi^4(n) = (\alpha_n + i\beta_n)^4 = (1 - 8\alpha_n^2 \beta_n^2) + 4i\alpha_n \beta_n (\alpha_n^2 - \beta_n^2).$$

Obviously,

$$-1 \leq 1 - 8\alpha_n^2 \beta_n^2 \leq 1,$$

$$-1 \leq 4\alpha_n \beta_n (\alpha_n^2 - \beta_n^2) \leq 1.$$

So, we have

$$\sqrt{\sum_{n=1}^{\infty} \frac{\chi^4(n)}{n^2}} \cdot \sqrt{\sum_{n=1}^{\infty} \frac{1}{n^{4it}}} \cdot \cos(\widehat{\vec{N}_3, \vec{N}_4}) =$$

$$= \sqrt{\sum_{n=1}^{\infty} \frac{1 - 8\alpha_n^2 \beta_n^2}{n^2} + i \cdot \sum_{n=1}^{\infty} \frac{4\alpha_n \beta_n (\alpha_n^2 - \beta_n^2)}{n^2}} \cdot \sqrt{\sum_{n=1}^{\infty} \frac{1}{n^{4it}}} \cdot \cos(\widehat{\vec{N}_3, \vec{N}_4}) = 0.$$

Take $(\widehat{\vec{N}_3, \vec{N}_4}) = \frac{\pi}{2}$, we have

$$\sqrt{\sum_{n=1}^{\infty} \frac{1 - 8\alpha_n^2 \beta_n^2}{n^2} + i \cdot \sum_{n=1}^{\infty} \frac{4\alpha_n \beta_n (\alpha_n^2 - \beta_n^2)}{n^2}} \cdot \sqrt{\sum_{n=1}^{\infty} \frac{1}{n^{4it}}} \cdot 0 = 0. \qquad (10)$$

Because

$$-\frac{\pi^2}{6} = -\sum_{n=1}^{\infty}\frac{1}{n^2} < \sum_{n=1}^{\infty}\frac{1-8\alpha_n^2\beta_n^2}{n^2} < \sum_{n=1}^{\infty}\frac{1}{n^2} = \frac{\pi^2}{6},$$

$$-\frac{\pi^2}{6} = -\sum_{n=1}^{\infty}\frac{1}{n^2} < \sum_{n=1}^{\infty}\frac{4\alpha_n\beta_n(\alpha_n^2-\beta_n^2)}{n^2} < \sum_{n=1}^{\infty}\frac{1}{n^2} = \frac{\pi^2}{6}$$

therefore, the complex series

$$\sum_{n=1}^{\infty}\frac{1-8\alpha_n^2\beta_n^2}{n^2} + i\cdot\sum_{n=1}^{\infty}\frac{4\alpha_n\beta_n(\alpha_n^2-\beta_n^2)}{n^2}$$

is convergent.

In above expression (10), when $t = 0$, although the series

$$\sum_{n=1}^{\infty}\sin(4t\ln n) = \sum_{n=1}^{\infty}\sin(4\cdot 0\cdot\ln n)$$
$$= 0+0+0+\cdots+0+\cdots = 0,$$

but because the series

$$\sum_{n=1}^{\infty}\cos(4t\ln n) = \sum_{n=1}^{\infty}\cos(4\cdot 0\cdot\ln n)$$
$$= 1+1+1+\cdots+1+\cdots = \infty,$$

therefore, the complex series

$$\frac{1}{1^{4it}}+\frac{1}{2^{4it}}+\frac{1}{3^{4it}}+\cdots+\frac{1}{n^{4it}}+\cdots$$
$$= \sum_{n=1}^{\infty}[\cos(4t\ln n) - i\sin(4t\ln n)]$$
$$= (1+1+1+\cdots+1+\cdots) -$$
$$- i\cdot(0+0+0+\cdots+0+\cdots) = \infty.$$

Because "infinite" multiplied by zero doesn't equal to zero, therefore

$$\sqrt{\sum_{n=1}^{\infty}\frac{1-8\alpha_n^2\beta_n^2}{n^2}+i\cdot\sum_{n=1}^{\infty}\frac{4\alpha_n\beta_n(\alpha_n^2-\beta_n^2)}{n^2}}\cdot\sqrt{\sum_{n=1}^{\infty}\frac{1}{n^{4it}}}\cdot 0$$

$$=\sqrt{\sum_{n=1}^{\infty}\frac{1-8\alpha_n^2\beta_n^2}{n^2}+i\cdot\sum_{n=1}^{\infty}\frac{4\alpha_n\beta_n(\alpha_n^2-\beta_n^2)}{n^2}}\cdot\sqrt{\sum_{n=1}^{\infty}[\cos(4t\ln n)-i\sin(4t\ln n)]}\cdot 0$$

$$=\sqrt{\sum_{n=1}^{\infty}\frac{1-8\alpha_n^2\beta_n^2}{n^2}+i\cdot\sum_{n=1}^{\infty}\frac{4\alpha_n\beta_n(\alpha_n^2-\beta_n^2)}{n^2}}\cdot\sqrt{\infty}\cdot 0\neq 0.$$

So, $\text{Re}(s)=\frac{1}{2}$ $(t=0)$ is not a non-trivial zero of the Dirichlet $L$ function $L(s,\chi)$, it must be removed.

When $t\neq 0$, although the series $\sum_{n=1}^{\infty}\cos(4t\ln n)$ and series $\sum_{n=1}^{\infty}\sin(4t\ln n)$ are all divergent, but they are all **oscillatory** and **bounded**, therefore we can know the complex series

$$\sum_{n=1}^{\infty}\frac{1}{n^{4it}}=\sum_{n=1}^{\infty}[\cos(4t\ln n)-i\sin(4t\ln n)]$$

is also **oscillatory** and **bounded**.

From following figure:

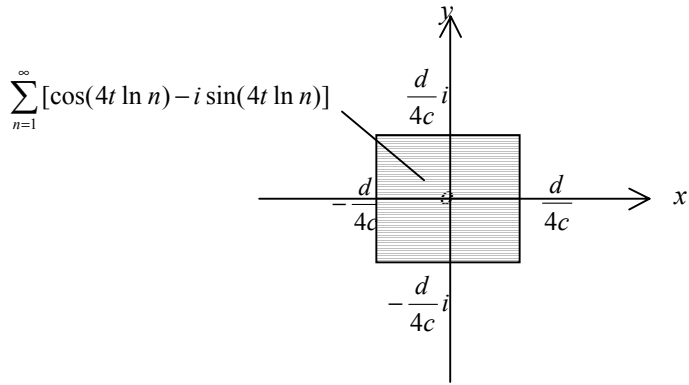

Fig. 4

we can know the change scope of the complex series $\sum_{n=1}^{\infty}[\cos(4t\ln n)-i\sin(4t\ln n)]$ is the square.

Obviously, the finite complex numbers multiplied by zero equal to zero, therefore we have

$$\sqrt{\sum_{n=1}^{\infty}\frac{1-8\alpha_n^2\beta_n^2}{n^2}+i\cdot\sum_{n=1}^{\infty}\frac{4\alpha_n\beta_n(\alpha_n^2-\beta_n^2)}{n^2}}\cdot\sqrt{\sum_{n=1}^{\infty}\frac{1}{n^{4it}}}\cdot 0$$

$$=\sqrt{\sum_{n=1}^{\infty}\frac{1-8\alpha_n^2\beta_n^2}{n^2}+i\cdot\sum_{n=1}^{\infty}\frac{4\alpha_n\beta_n(\alpha_n^2-\beta_n^2)}{n^2}}\cdot\sqrt{\sum_{n=1}^{\infty}[\cos(4t\ln n)-i\sin(4t\ln n)]}\cdot 0$$

$$=0, \qquad\qquad\qquad (t\neq 0)$$

**This explains that the imaginative part $t$ of the non-trivial zeroes of Dirichlet $L$ function takes all values of the straight line $\mathrm{Re}(s)=\dfrac{1}{2}$ except $t=0$, namely,**

$$t\in(-\infty,\ 0)\cup(0,\ +\infty).$$

On the other hand, from the expression (6), we have

$$L(s,\chi)=\sum_{n=1}^{\infty}\frac{\chi^2(n)}{n^{2s}}=(\frac{1}{1^{2\sigma}},\ \frac{1}{2^{2\sigma}},\ \frac{1}{3^{2\sigma}},\ \cdots,\ \frac{1}{n^{2\sigma}},\ \cdots)\cdot(\frac{\chi^2(1)}{1^{2it}},\ \frac{\chi^2(2)}{2^{2it}},\ \frac{\chi^2(3)}{3^{2it}},\ \cdots,\ \frac{\chi^2(n)}{n^{2it}},\ \cdots)$$

$$=\sqrt{\sum_{n=1}^{\infty}\frac{1}{n^{4\sigma}}}\cdot\sqrt{\sum_{n=1}^{\infty}\frac{\chi^4(n)}{n^{4it}}}\cdot\cos(\widehat{\vec{N}_3,\vec{N}_4})$$

$$=\sqrt{\sum_{n=1}^{\infty}\frac{1}{n^{4\sigma}}}\cdot\sqrt{\sum_{n=1}^{\infty}\chi^4(n)\cdot e^{-4it\log n}}\cdot 0$$

$$=\sqrt{\sum_{n=1}^{\infty}\frac{1}{n^{4\sigma}}}\cdot\sqrt{\sum_{n=1}^{\infty}\chi^4(n)(\cos(4t\log n)-i\sin(4t\log n))}\cdot 0$$

$$=\sqrt{\sum_{n=1}^{\infty}\frac{1}{n^{4\sigma}}}\cdot\sqrt{\sum_{n=1}^{\infty}((1-8\alpha_n^2\beta_n^2)+i(4\alpha_n\beta_n(\alpha_n^2-\beta_n^2)))(\cos(4t\log n)-i\sin(4t\log n))}\cdot 0$$

$$=\sqrt{\sum_{n=1}^{\infty}\frac{1}{n^{4\sigma}}}\cdot\sqrt{\sum_{n=1}^{\infty}\begin{pmatrix}((1-8\alpha_n^2\beta_n^2)\cos(4t\log n)+4\alpha_n\beta_n(\alpha_n^2-\beta_n^2)\sin(4t\log n))+\\ +i(4\alpha_n\beta_n(\alpha_n^2-\beta_n^2)\cos(4t\log n)-(1-8\alpha_n^2\beta_n^2)\sin(4t\log n))\end{pmatrix}}\cdot 0=0.$$

Obviously, when $t=0$, the complex series

$$\sum_{n=1}^{\infty}\begin{pmatrix}((1-8\alpha_n^2\beta_n^2)\cos(4t\log n)+4\alpha_n\beta_n(\alpha_n^2-\beta_n^2)\sin(4t\log n))+\\ +i(4\alpha_n\beta_n(\alpha_n^2-\beta_n^2)\cos(4t\log n)-(1-8\alpha_n^2\beta_n^2)\sin(4t\log n))\end{pmatrix}=\infty,$$

because "infinite" multiplied by zero doesn't equal to zero, therefore

$$\sqrt{\sum_{n=1}^{\infty}\frac{1}{n^{4\sigma}}} \cdot \sqrt{\sum_{n=1}^{\infty}\begin{pmatrix}((1-8\alpha_n^2\beta_n^2)\cos(4t\log n)+4\alpha_n\beta_n(\alpha_n^2-\beta_n^2)\sin(4t\log n))+\\+i(4\alpha_n\beta_n(\alpha_n^2-\beta_n^2)\cos(4t\log n)-(1-8\alpha_n^2\beta_n^2)\sin(4t\log n))\end{pmatrix}} \cdot \cos(\widehat{\vec{N}_3,\vec{N}_4})$$

$$=\sqrt{\sum_{n=1}^{\infty}\frac{1}{n^{4\sigma}}} \cdot \sqrt{\infty} \cdot 0 \neq 0.$$

So, $\mathrm{Re}(s)=\frac{1}{2}$ $(t=0)$ is not a non-trivial zero of the Dirichlet $L$ function $L(s,\chi)$, it must be removed.

When $t \neq 0$, the complex series

$$\sum_{n=1}^{\infty}\begin{pmatrix}((1-8\alpha_n^2\beta_n^2)\cos(4t\log n)+4\alpha_n\beta_n(\alpha_n^2-\beta_n^2)\sin(4t\log n))+\\+i(4\alpha_n\beta_n(\alpha_n^2-\beta_n^2)\cos(4t\log n)-(1-8\alpha_n^2\beta_n^2)\sin(4t\log n))\end{pmatrix}$$

is also **oscillatory** and **bounded**.

Obviously, the finite complex numbers multiplied by zero equal to zero, therefore, we have

$$\sqrt{\sum_{n=1}^{\infty}\frac{1}{n^{4\sigma}}} \cdot \sqrt{\sum_{n=1}^{\infty}\begin{pmatrix}((1-8\alpha_n^2\beta_n^2)\cos(4t\log n)+4\alpha_n\beta_n(\alpha_n^2-\beta_n^2)\sin(4t\log n))+\\+i(4\alpha_n\beta_n(\alpha_n^2-\beta_n^2)\cos(4t\log n)-(1-8\alpha_n^2\beta_n^2)\sin(4t\log n))\end{pmatrix}} \cdot 0 = 0.$$

**This explains that the imaginative part *t* of the non-trivial zeroes of Dirichlet *L* function takes all values of the straight line** $\mathrm{Re}(s)=\frac{1}{2}$ **except** $t=0$, namely

$$t \in (-\infty,\ 0) \cup (0,\ +\infty).$$

# Appendix

## The inner product formula between two vectors in real space can be extended formally to complex space

**Dear Mr. Referee,**

**Thank you for your review.**

**First, please observe following examples.**

**Example 1**  $A\{1+i,\ 3\},\ B\{-i,\ 2i\},\ C\{1,\ -i\}$.

$$\overrightarrow{AB} = \{-1-2i,\ -3+2i\},\ \overrightarrow{AC} = \{-i,\ -3-i\},\ \overrightarrow{BC} = \{1+i,\ -3i\}.$$

$$\left|\overrightarrow{AB}\right| = \sqrt{(-1-2i)^2 + (-3+2i)^2} = \sqrt{2-8i},$$

$$\left|\overrightarrow{AC}\right| = \sqrt{(-i)^2 + (-3-i)^2} = \sqrt{7+6i},$$

$$\left|\overrightarrow{BC}\right| = \sqrt{(1+i)^2 + (-3i)^2} = \sqrt{-9+2i}.$$

According to **Cosine theorem**

$$\left|\overrightarrow{AB}\right|^2 + \left|\overrightarrow{AC}\right|^2 - \left|\overrightarrow{BC}\right|^2 = 2\left|\overrightarrow{AB}\right|\left|\overrightarrow{AC}\right|\cos(\overrightarrow{AB}, \overrightarrow{AC}),$$

we have

$$\frac{\left|\overrightarrow{AB}\right|^2 + \left|\overrightarrow{AC}\right|^2 - \left|\overrightarrow{BC}\right|^2}{2} = \frac{(2-8i)+(7+6i)-(-9+2i)}{2} = 9-2i.$$

On the other hand, we have

$$\overrightarrow{AB} \cdot \overrightarrow{AC} = \left|\overrightarrow{AB}\right|\left|\overrightarrow{AC}\right|\cos(\overrightarrow{AB}, \overrightarrow{AC}) = (-1-2i)(-i)+(-3+2i)(-3-i) = 9-2i.$$

According to the **Cosine theorem**

$$\left|\overrightarrow{AC}\right|^2 + \left|\overrightarrow{BC}\right|^2 - \left|\overrightarrow{AB}\right|^2 = 2\left|\overrightarrow{AC}\right|\left|\overrightarrow{BC}\right|\cos(\overrightarrow{AC}, \overrightarrow{BC})$$

$$\frac{\left|\overrightarrow{AC}\right|^2 + \left|\overrightarrow{BC}\right|^2 - \left|\overrightarrow{AB}\right|^2}{2} = \frac{(7+6i)+(-9+2i)-(2-8i)}{2} = -2+8i.$$

On the other hand, we have

$$\overrightarrow{AC} \cdot \overrightarrow{BC} = \left|\overrightarrow{AB}\right|\left|\overrightarrow{AC}\right|\cos(\overrightarrow{AB}, \overrightarrow{AC}) = (-i)(1+i)+(-3-i)(-3i) = -2+8i.$$

$$S_{\triangle ABC} = \frac{1}{2}\left|\overrightarrow{AB}\right|\left|\overrightarrow{AC}\right|\sin\theta_1 = \frac{1}{2}\sqrt{2-8i}\cdot\sqrt{7+6i}\cdot\frac{\sqrt{-15-8i}}{\sqrt{62-44i}} = \frac{1}{2}\sqrt{-15-8i},$$

$$S_{\triangle ACB} = \frac{1}{2}\left|\overrightarrow{AC}\right|\left|\overrightarrow{BC}\right|\sin\theta_{21} = \frac{1}{2}\sqrt{7+6i}\cdot\sqrt{-9+2i}\cdot\frac{\sqrt{-15-8i}}{\sqrt{-75-40i}} = \frac{1}{2}\sqrt{-15-8i}.$$

**Example 2** $A\{1+i, 1-i, 2i\}$, $B\{1-i, 1+i, -2i\}$, $C\{1, 0, i\}$.

$\overrightarrow{AB} = \{-2i, 2i, -4i\}$, $\overrightarrow{AC} = \{-i, -1+i, -i\}$, $\overrightarrow{BC} = \{i, -1-i, 3i\}$.

$\left|\overrightarrow{AB}\right| = \sqrt{(-2i)^2 + (2i)^2 + (-4i)^2} = \sqrt{-24},$

$\left|\overrightarrow{AC}\right| = \sqrt{(-i)^2 + (-1+i)^2 + (-i)^2} = \sqrt{-2-2i},$

$\left|\overrightarrow{BC}\right| = \sqrt{(i)^2 + (-1-i)^2 + (3i)^2} = \sqrt{-10+2i}.$

According to the **Cosine theorem**

$$\left|\overrightarrow{AB}\right|^2 + \left|\overrightarrow{AC}\right|^2 - \left|\overrightarrow{BC}\right|^2 = 2\left|\overrightarrow{AB}\right|\left|\overrightarrow{AC}\right|\cos(\overrightarrow{AB}, \overrightarrow{AC}),$$

we have

$$\frac{|\overrightarrow{AB}|^2 + |\overrightarrow{AC}|^2 - |\overrightarrow{BC}|^2}{2} = \frac{(-24) + (-2 - 2i) - (-10 + 2i)}{2} = -8 - 2i.$$

On the other hand, we have

$$\overrightarrow{AB} \cdot \overrightarrow{AC} = |\overrightarrow{AB}||\overrightarrow{AC}|\cos(\overrightarrow{AB}, \overrightarrow{AC}) = (-2i)(-i) + (2i)(-1+i) + (-4i)(-i) = -8 - 2i.$$

According to the **Cosine theorem**

$$|\overrightarrow{AC}|^2 + |\overrightarrow{BC}|^2 - |\overrightarrow{AC}|^2 = 2|\overrightarrow{AC}||\overrightarrow{BC}|\cos(\overrightarrow{AC}, \overrightarrow{BC}),$$

we have

$$\frac{|\overrightarrow{AC}|^2 + |\overrightarrow{BC}|^2 - |\overrightarrow{AB}|^2}{2} = \frac{(-2 - 2i) + (-10 + 2i) - (-24)}{2} = 6.$$

On the other hand, we have

$$\overrightarrow{AC} \cdot \overrightarrow{BC} = |\overrightarrow{AC}||\overrightarrow{BC}|\cos(\overrightarrow{AC}, \overrightarrow{BC}) = (-i)(i) + (-1+i)(-1-i) + (-i)(3i) = 6.$$

$$S_{\triangle ABC} = \frac{1}{2}|\overrightarrow{AB}||\overrightarrow{AC}|\sin\theta_1 = \frac{1}{2}\sqrt{-24} \cdot \sqrt{-2-2i} \cdot \frac{\sqrt{-12+16i}}{\sqrt{48+48i}} = \sqrt{-3+4i},$$

$$S_{\triangle ACB} = \frac{1}{2}|\overrightarrow{AC}||\overrightarrow{BC}|\sin\theta_{21} = \frac{1}{2}\sqrt{-2-2i} \cdot \sqrt{-10+2i} \cdot \frac{\sqrt{-12+16i}}{\sqrt{24+16i}} = \sqrt{-3+4i}.$$

**Example 3** $A\{8i,\ 14,\ 8-i,\ 1\}$, $B\{6,\ 15i,\ 17,\ -8\}$, $C\{3-i,\ 10+7i,\ 11,\ 3i\}$.

$\overrightarrow{AB} = \{6 - 8i,\ -14 + 15i,\ 9 + i,\ -9\}$, $\overrightarrow{AC} = \{3 - 9i,\ -4 + 7i,\ 3 + i,\ -1 + 3i\}$,

$\overrightarrow{BC} = \{-3 - i,\ 10 - 8i,\ -6,\ 8 + 3i\}$.

$|\overrightarrow{AB}| = \sqrt{(6-8i)^2 + (-14+5i)^2 + (9+i)^2 + (-9)^2} = \sqrt{104 - 498i}$,

$|\overrightarrow{AC}| = \sqrt{(3-9i)^2 + (-4+7i)^2 + (3+i)^2 + (-1+3i)^2} = \sqrt{-105 - 110i}$,

$|\overrightarrow{BC}| = \sqrt{(-3-i)^2 + (10-8i)^2 + (-6)^2 + (8+3i)^2} = \sqrt{135 - 106i}$.

According to the **Cosine theorem**

$$\left|\overrightarrow{AB}\right|^2 + \left|\overrightarrow{AC}\right|^2 - \left|\overrightarrow{BC}\right|^2 = 2\left|\overrightarrow{AB}\right|\left|\overrightarrow{AC}\right|\cos(\overrightarrow{AB},\overrightarrow{AC}),$$

we have

$$\frac{\left|\overrightarrow{AB}\right|^2 + \left|\overrightarrow{AC}\right|^2 - \left|\overrightarrow{BC}\right|^2}{2} = \frac{(104-498i)+(-105-110i)-(-135-106i)}{2} = -68-251i.$$

On the other hand, we have

$$\overrightarrow{AB}\cdot\overrightarrow{AC} = \left|\overrightarrow{AB}\right|\left|\overrightarrow{AC}\right|\cos(\overrightarrow{AB},\overrightarrow{AC}) =$$
$$= (6-8i)(3-9i)+(-14+15i)(-4+7i)+(9+i)(3+i)+(-9)(-1+3i) = -68-251i.$$

According to the **Cosine theorem**

$$\left|\overrightarrow{AC}\right|^2 + \left|\overrightarrow{BC}\right|^2 - \left|\overrightarrow{AC}\right|^2 = 2\left|\overrightarrow{AC}\right|\left|\overrightarrow{BC}\right|\cos(\overrightarrow{AC},\overrightarrow{BC}),$$

we have

$$\frac{\left|\overrightarrow{AC}\right|^2 + \left|\overrightarrow{BC}\right|^2 - \left|\overrightarrow{AB}\right|^2}{2} = \frac{(-105-110i)+(135-106i)-(104-498i)}{2} = -37+141i.$$

On the other hand, we hand

$$\overrightarrow{AC}\cdot\overrightarrow{BC} = \left|\overrightarrow{AC}\right|\left|\overrightarrow{BC}\right|\cos(\overrightarrow{AC},\overrightarrow{BC}) =$$
$$= (3-9i)(-3-i)+(-4+7i)(10-8i)+(-6)(3+i)+(-1+3i)(8+3i) = -37+141i.$$

$$S_{\Delta ABC} = \frac{1}{2}\left|\overrightarrow{AB}\right|\left|\overrightarrow{AC}\right|\sin\theta_1 = \frac{1}{2}\sqrt{104-498i}\cdot\sqrt{-105-110i}\cdot\frac{\sqrt{-7323+6714i}}{\sqrt{-65700+40850i}}$$
$$= \frac{1}{2}\sqrt{-7323+6714i},$$

$$S_{\Delta ACB} = \frac{1}{2}\left|\overrightarrow{AC}\right|\left|\overrightarrow{BC}\right|\sin\theta_{21} = \frac{1}{2}\sqrt{-105-110i}\cdot\sqrt{135-106i}\cdot\frac{\sqrt{-7323+6714i}}{\sqrt{-25835-3720i}}$$
$$= \frac{1}{2}\sqrt{-7323+6714i}.$$